\documentclass[12pt]{iopart}
\pdfoutput=1 
\usepackage{graphicx}
\begin{document}

\title{Continuum Cascade Model of Directed Random Graphs: Traveling Wave Analysis}

\author{Yoshiaki  Itoh}
\address{The Institute of Statistical Mathematics, 
10-3 Midori-cho, Tachikawa, Tokyo 190-8562, Japan}
\author{P.~L.~Krapivsky}
\address{Department of Physics, Boston University, Boston, MA 02215, USA}

\begin{abstract}
We study a class of directed random graphs. In these graphs, the interval $[0,x]$ is the vertex set, and from each $y\in [0,x]$, directed links are drawn to points in the interval $(y,x]$ which are chosen uniformly with density one. We analyze the length of the longest directed path starting from the origin. In the $x\to\infty$ limit, we employ traveling wave techniques to extract the asymptotic behavior of this quantity.  We also study the size of a cascade tree composed of vertices which can be reached via directed paths starting at the origin. 
\end{abstract}


\maketitle

\section{Introduction}

A random graph is a set of vertices that are connected by random links \cite{sr,er,bol,s,dk,jklp,jlr}. 
Random graphs underlie numerous natural phenomena ranging from polymerization \cite{pf,book} to the
spread of infectious diseases \cite{mej}, and they also have applications to transportation systems, 
electrical distribution systems, the Internet, the world-wide web, social networks, etc. \cite{gc,DM03,mejn}.  

In random graph models, links are usually treated as undirected. In a growing number of applications, however, directionality plays a prominent role. One example is modeling of the web growth \cite{broder,krr,kr,DM03,mejn}. In modeling of food webs directionality (reflecting predation) is even more crucial. Food webs are directed graphs with vertexes $\{0,1,\ldots,m\}$ labeling different species. The presence of the directed link $(i, j)$ indicates that species $i$ is eaten by species $j$. Usually in food web only links $(i, j)$ with $i<j$ are allowed. (Loops $(i,i)$ which would account cannibalism are ignored; the directed link $(i, j)$ with $i>j$ could e.g. represent predation on the young of the `stronger' species $i$ by adults of species $j$, but such links are also disregarded in most  models.) The simplest cascade model  \cite{c,cbn,cn,cbn2,cn2,n} generates a food web at random, namely for each pair of species $i$ and $j$ with $i<j$ the directed link  $(i, j)$ is drawn at random with a certain predation probability $c$. A number of questions, particularly those related to the maximal length of food chains, have been investigated in the framework of the this cascade model.  For instance, what is the length (the number of links) of the longest direct path starting from the {\em basal} species (vertex 0)? A dual question concerns the length of the longest path finishing at the {\em top} species (vertex $m$). One can also ask about the length $M$ of the longest path irrespectively on the first and last species. 

The simplest cascade model is a kind of `standard model' in the subject, and it had been widely used to interpret ecological data on community food webs \cite{cbn}. The standard cascade model provides a very natural mechanism for generating directed random graphs and the same model has been suggested in other contexts, e.g. as a model of parallel computation  \cite{gel,newman_94} in which the presence of the directed link $(i, j)$ with $i<j$ indicates that task $i$ must be performed before task $j$. For a parallel computation in which each task takes a unit of time, the processing time will be equal $M+1$ (where $M$ is the length of the the longest path). 

Food webs typically involve a huge number of species\footnote{Small food webs tend to reflect our ignorance rather than reality.}, while the average predation per species is usually not too large. Hence it is interesting to investigate large food webs with small predation probability, more precisely the scaling limit 
\begin{equation}
\label{scaling}
m\to\infty, \quad c\to 0, \quad cm= x = {\rm finite}
\end{equation}
This suggests to study a continuum cascade model where the vertex set is the interval $[0,x]$. For each species $y$, the number of predator species is random, such species are chosen at random from the interval $(y,x]$  according to the Poisson distribution with unit density. The Poisson distribution immediately follows from the binomial distribution (characterizing the discrete cascade model) in the scaling limit \eref{scaling}. This cascade model is the minimalist continuum model of directed random graphs. Simple models tend to arise in various unrelated subjects and they are interesting on purely intellectual grounds. Nevertheless, for concreteness in the following exposition we shall often use the language of food webs. 

The rest of this article is organized as follows. In Sec.~\ref{model} we define the model, discuss its simplest properties, and derive a recurrence for the longest directed path starting from the origin. The asymptotic behavior of the solution to that recurrence is analyzed in the following sections \ref{elem} and \ref{t_wave}. In section \ref{size_tree} we discuss the total number of vertices in a cascade tree with the root at the origin; on the language of food webs it counts the basal species and species feeding on it, both directly and indirectly. 

\section{Continuum Cascade Model} 
\label{model}

The vertex set of our random graph is the interval $[0,x]$. In the illustrative picture below we draw only the vertex set and links from the cascade subgraph initiating at the origin (the open circle on the picture). Namely, we draw all links emanating from the origin indicating direct predation on the basal species (there are 3 such predators in the picture); then we draw all the links from these direct predators (4 such predators in the picture); etc. Links are drawn in a cascade manner thereby explaining the name of the model. 

\begin{picture}(100,120)(20,10)
\put(20,50){\line (420,0) {420}}
\thicklines
\put(100,50){\oval (160,50) [t]}
\put(150,50){\oval (260,80) [t]}
\put(200,50){\oval (360,120) [t]}
\put(210,50){\oval (60,20) [b]}
\put(260,50){\oval (160,50) [b]}
\put(300,50){\oval (240,70) [b]}
\put(275,50){\oval (70,25) [t]}
\put(300,50){\oval (120,40) [t]}
\put(340,50){\oval (120,40) [b]}
\put(375,50){\oval (70,25) [t]}
\put(20,50){\circle{5}}
\put(310,50){\circle*{5}}
\put(360,50){\circle*{5}}
\put(380,50){\circle*{5}}
\put(400,50){\circle*{5}}
\put(410,50){\circle*{5}}
\put(420,50){\circle*{5}}
\put(435,35){{\bf x}}
\end{picture}

Overall, in the above illustrative picture the cascade subgraph is a tree with 10 links and 11 vertices. Six of these vertices (closed circles on the picture) are terminal, that is, there are no links emanating from them. 
Every cascade subgraph is a tree; the size and the number of terminal vertices in cascade trees fluctuate from realization to realization. 

Terminal vertices represent top predators on the language of food webs. It is easy to compute the fraction of top predators:
\begin{equation}
T = \frac{1}{x}\int_0^x dy\,e^{-y} = \frac{1-e^{-x}}{x}
\end{equation}
The fraction of bottom preys\footnote{Bottom preys are often called basal species. We reserve the term `basal species' only for the species at the origin which, according to the definition of the continuum cascade model, can never be a predator independently on the choice of links.}, that is, species who do not eat other species, is the same. The overlap of the sets of top predators and bottom preys (one can call them neutral species) is non-empty, the fraction of neutral species is
\begin{equation}
N = \frac{1}{x}\int_0^x dy\,e^{-y}\,e^{-(x-y)} = e^{-x}
\end{equation}

\begin{figure}
\hspace*{-0.1in}
  \centering
\includegraphics[scale=0.45]{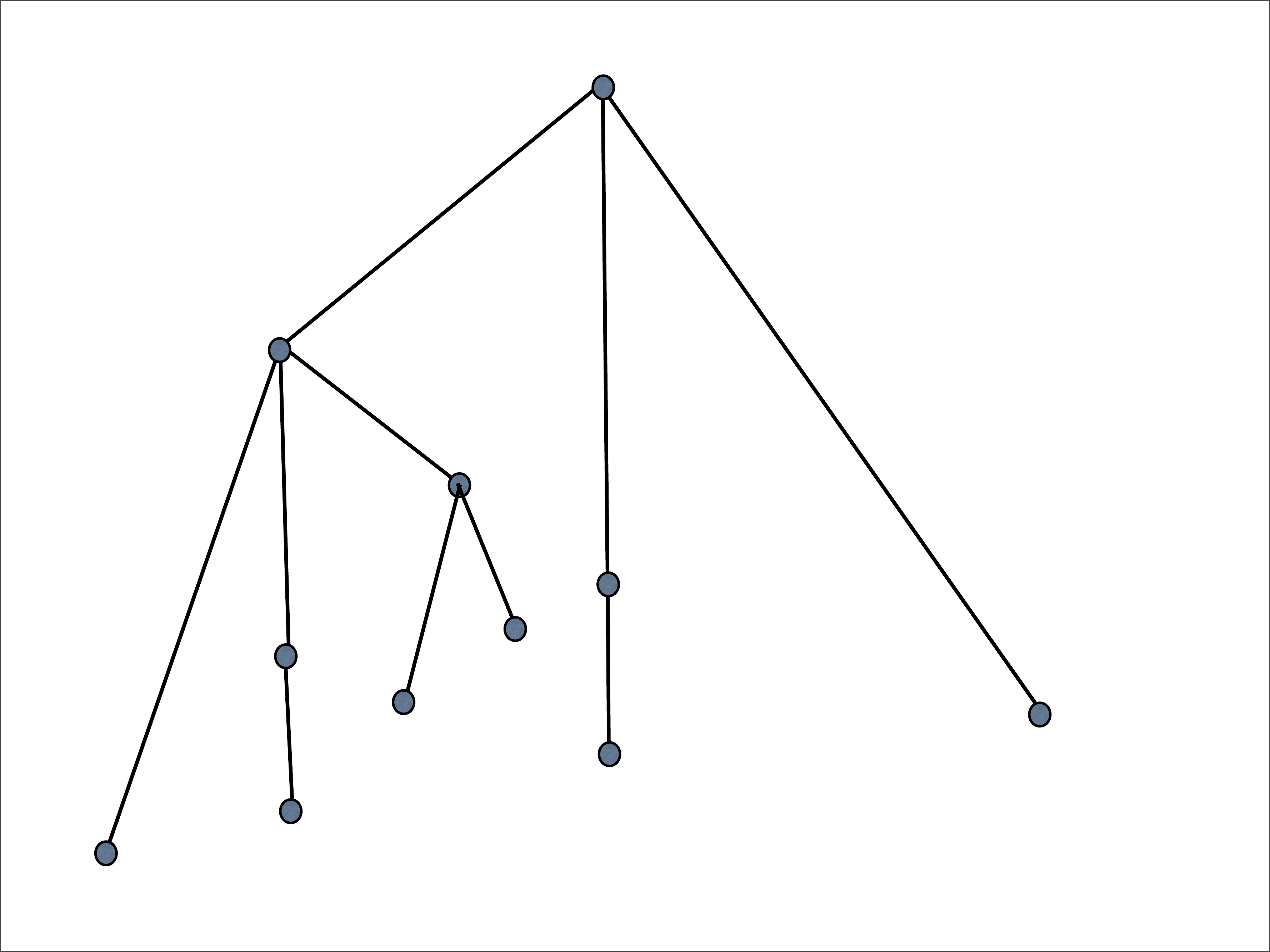}
\caption{The cascade tree with the basal species (the vertex at the top) playing the role of the root. The height of this cascade tree is equal to 3.} 
\label{fig:food_web}
\end{figure}

We now turn to more subtle properties of the continuum cascade model which are related to the cascade tree. This tree is finite and it varies from realization to realization; accordingly, the properties of the cascade tree are probabilistic. To define these properties it is convenient to utilize a more traditional way of plotting trees; the cascade tree pictured above is presented on Fig.~\ref{fig:food_web}. This figure resembles binary search trees and both the relevant properties of binary search trees and the methods used in analysis of binary search trees \cite{Drmota,si,ho,jn,di,r,fo,d,km,bkm,mk,sz,itoh} are useful in our situation. For instance, the height of the binary search trees has attracted a lot of attention, and a traveling wave analysis \cite{km,bkm,mk} has provided a very efficient way of tackling the asymptotic (in number of vertices of the tree) behavior of the height. In the present problem, the height is indeed an interesting quantity, namely it is the length of longest chain from the basal species to the bottom of the cascade tree, and the traveling wave analysis will be helpful as well.

We now establish a recurrence relation for the height distribution. The height is a non-negative integer. It is convenient to work with the cumulative distribution 
\begin{equation}
P_n(x) = {\rm Prob}({\rm height} \leq n)
\end{equation}

The basal species is the terminal vertex with probability $e^{-x}$, and therefore
\begin{equation}
\label{P0}
P_0(x)=e^{-x}
\end{equation}
For $n\geq 1$, 
\begin{eqnarray*}
P_n(x)&=&  
\sum _{k\geq 0} \frac{x^k }{k!}\,e^{-x}  \int _0^x \cdots  
\int _0^x  P_{n-1}(y_1)\cdots P_{n-1}(y_k)\,\frac{dy_1}{x}\cdots \frac{dy_k}{x}\\
&=&  e^{-x}  \sum _{k\geq 0} \frac{1}{k!} \left(   \int _0^x  P_{n-1}(y)\,dy \right)^k,
\end{eqnarray*} 
where the first line accounts for any possible number $k\geq 0$ of links emanating from the origin and finishing at all possible points $x-y_j$. There is no `interaction' between different branches of the cascade tree, so we merely must assure that all cascade trees originating at $x-y_j$ have heights not exceeding $n-1$. Computing the sum in the above equation we arrive at our main recurrence
\begin{equation}
\label{eq:height}
P_n(x) = \exp\left[-x + \int _0^x  P_{n-1}(y)\,dy)\right]
\end{equation}

\begin{figure}
  \centering
  \includegraphics[scale=1.2]{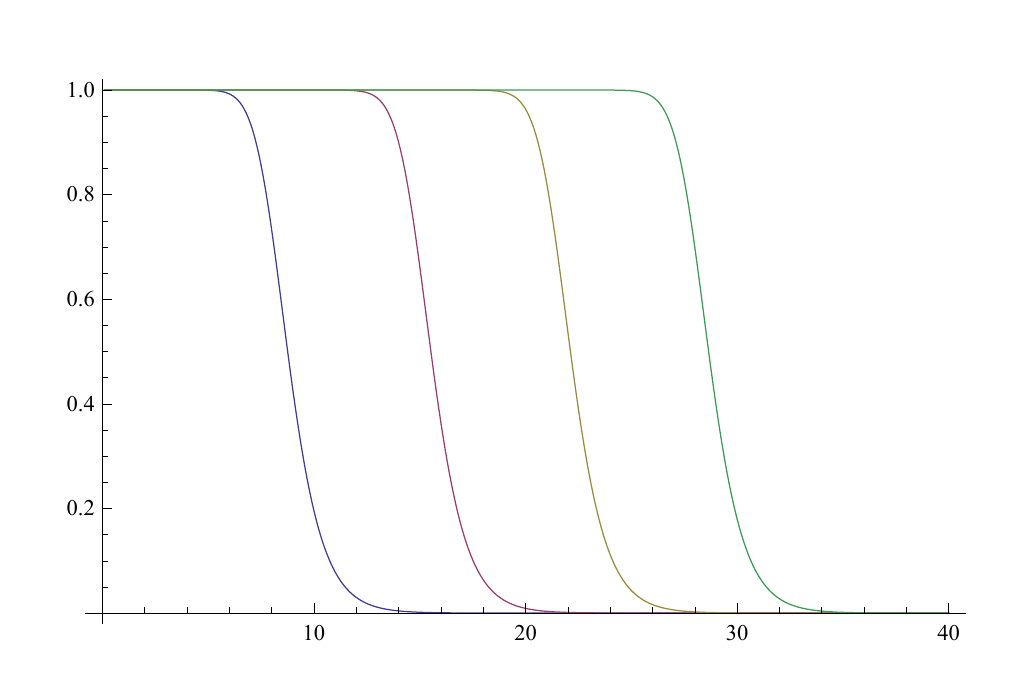}
  \caption{The distribution $P_n(x)$ versus $x$ obtained by iteration the recurrence 
  \eref{eq:height}. Shown is $P_n(x)$ for $n=20, 40, 60, 80$ (left to right). Iterations were performed 
  using {\it Mathematica}. The observed velocity of the traveling wave is in a good agreement with 
  the theoretical prediction \eref{front}.} 
\label{fig:TW}
\end{figure}

Starting with \eref{P0} we find
\begin{equation}
\label{P_1}   
P_1(x) = \exp\!\left[-x+1-e^{-x}\right].
\end{equation}
One can recursively determine $P_2$, then $P_3$; analytical expressions for  $P_n$ become very cumbersome as $n$ increases. Fortunately, in the large `time' limit, $n\gg 1$, the behavior greatly simplifies, namely the solution acquires a traveling wave form [see Fig.~\ref{fig:TW}],
\begin{equation}
\label{TW}
P_n(x)\to \Pi(\xi), \quad \xi=x-x_f,
\end{equation}
with the  front position growing linearly with `velocity' equal to $e^{-1}$:
\begin{equation}
\label{front}
x_f\simeq vn, \quad v=\frac{1}{e}
\end{equation}

The traveling wave profile $\Pi(\xi)$ decreases monotonically from 1 to 0 as $\xi$ increases 
from $-\infty$ to $\infty$. More precisely,  
\begin{equation}
\label{large}
    \Pi(\xi)\propto e^{-\xi}\quad {\rm when} \quad \xi\to\infty
\end{equation}
and
\begin{equation}
\label{small}
   1-\Pi(\xi)\propto e^{e\xi}\quad  {\rm when} \quad \xi\to -\infty
\end{equation}
In the next section we give an elementary argument which allows one to understand \eref{front}. A more comprehensive traveling wave analysis that leads to above results  is presented in section \ref{t_wave}.

\section{Elementary derivation of the traveling wave velocity}
\label{elem}

Let us begin with the behavior of $P_n(x)$ for small $x$. Expanding $P_0(x)$ and $P_1(x)$, see 
equation  \eref{P_1}, we obtain 
\begin{eqnarray*}
P_0 = 1 - x + \frac{1}{2}x^2 - \frac{1}{6}x^3 + \ldots, \quad 
P_1 = 1 - \frac{1}{2}x^2 + \frac{1}{6}x^3 + \frac{1}{12}x^4 + \ldots
\end{eqnarray*}
Using equation \eref{eq:height} we recurrently determine the expansions of the following $P_n$ to yield
\begin{equation}
\label{Pn_2}
P_n = 1 - \frac{1}{(n+1)!}x^{n+1} + \ldots
\end{equation}
This result is easy to prove by induction. One can continue this expansion, e.g.,
\begin{equation}
\label{Pn_4}
P_n = 1 - \frac{1}{(n+1)!}x^{n+1} + \frac{1}{(n+2)!}x^{n+2} + \frac{2}{(n+3)!}x^{n+3} + \ldots
\end{equation}
is valid for all $n\geq 1$; this is also easily proven by induction. 

Let us now estimate the front position $x_f$ from the criterion $P_n(x_f)=\frac{1}{2}$. Keeping only two terms as in equation \eref{Pn_2} we obtain $x_f^{n+1}=\frac{1}{2}(n+1)!$, or $x_f=\frac{n+1}{e}$ in the leading order. What will happen if we keep e.g. four terms in the expansion? Using the same criterion $P_n(x_f)=\frac{1}{2}$ in conjunction with equation \eref{Pn_4} we get
\begin{eqnarray*}
\frac{x_f^{n+1}}{(n+1)!}\left[1-\frac{x_f}{n+2}-\frac{2 x_f^2}{(n+2)(n+3)}\right]=\frac{1}{2}.
\end{eqnarray*}
In the leading order we recover the previous prediction $x_f=\frac{n+1}{e}$.  This does not prove equation \eref{front},  but at least it shows its consistency with the series \eref {Pn_4}.

\section{Traveling wave analysis: Velocity selection}
\label{t_wave} 

We want to understand the behavior of the recurrence
\begin{equation}
\label{rec_def}   
P_{n+1}(x) = \exp\!\left[-x+\int_0^x dy\,P_n(y)\right]
\end{equation}
when $x\gg 1$. We assume the convergence to traveling wave solution and the validity of the traveling wave ansatz  \eref{TW}. Numerical results strongly support this assumption  [see Fig.~\ref{fig:TW}] and show that the convergence is rather fast, that is, the asymptotic shape emerges already for not to large $n$. We also assume that $x_f=vn$ (for large $n$), but we do not specify $v$. The left-hand side of equation \eref{rec_def} becomes 
\begin{equation}
\label{LHS}
P_{n+1}(x) = \Pi(\xi-v)
\end{equation}
while the right-hand side of  equation \eref{rec_def} turns into 
\begin{equation}
\label{RHS_long}
\exp\!\left[-x+\int_{-nv}^\xi d\eta\,\Pi(\eta)\right], \qquad \eta=y-nv
\end{equation}
We will allow $|\xi|$ to be large, but we will always assume that $|\xi|\ll nv$. In this situation, 
the integral in equation \eref{RHS_long} can be simplified as follows:
\begin{eqnarray}
\label{RHS_int}
\int_{-nv}^\xi d\eta\,\Pi(\eta) &=& \int_{-nv}^0 d\eta\,\Pi(\eta) + \int_0^\xi d\eta\,\Pi(\eta) \nonumber\\
&=&  nv + \int_{-nv}^0 d\eta\,[\Pi(\eta)-1] + \int_0^\xi d\eta\,\Pi(\eta) \nonumber\\
&=&  nv + \int_0^\xi d\eta\,\Pi(\eta) - L + L_n, 
\end{eqnarray}
where we have used the shorthand notation
\begin{equation}
\label{L_def}
L = \int_{-\infty}^0 d\eta\,[1-\Pi(\eta)], \quad L_n = \int_{-\infty}^{-nv} d\eta\,[1-\Pi(\eta)]. 
\end{equation}
Since $\Pi(\eta)$ quickly approaches to 1 as $\eta\to -\infty$, we drop $L_n$ from \eref{RHS_int}; we shall justify this step {\em a posteriori}. Combining \eref{RHS_long} and \eref{RHS_int} we see that the right-hand side of \eref{rec_def} becomes
\begin{equation}
\label{RHS}
\exp\!\left[-\xi - L + \int_0^\xi d\eta\,\Pi(\eta)\right]. 
\end{equation}
Equating \eref{LHS}  and \eref{RHS} we arrive at 
\begin{equation}
\label{main}
\Pi(\xi-v) = \exp\!\left[-\xi - L + \int_0^\xi d\eta\,\Pi(\eta)\right].
\end{equation}
This is the governing equation for $\Pi(\xi)$. 

\subsection{Far ahead of the front: $\xi\to \infty$}

Since $\Pi(\xi)$ quickly approaches to zero as $\xi\to\infty$, equation \eref{main} gives 
$\Pi(\xi-v)\simeq e^{-\xi-L+R}$, where $R = \int_0^\infty d\eta\,\Pi(\eta)$. Thus we confirm \eref{large}; more precisely we get 
\begin{equation}
\label{large_asymp}
\Pi(\xi)\simeq e^{R-L-v}\, e^{-\xi} \quad {\rm when} \quad \xi\to\infty.
\end{equation} 

\subsection{Far behind the front: $\xi\to -\infty$}

It is more convenient to work with $\Phi(\xi) = 1-\Pi(\xi)$ rather than $\Pi(\xi)$. In terms of 
$\Phi(\xi)$,  equation \eref{main} becomes
\begin{equation}
\label{main_Phi}
1-\Phi(\xi-v) = \exp\!\left[ - L - \int_0^\xi d\eta\,\Phi(\eta)\right].
\end{equation}
From the definition \eref{L_def} we see that 
\begin{eqnarray*}
L = \int_{-\infty}^0 d\eta\,\Phi(\eta).
\end{eqnarray*}
Therefore we can re-write \eref{main_Phi} as
\begin{equation}
\label{chief}
1-\Phi(\xi-v) = \exp\!\left[ - \int_{-\infty}^\xi d\eta\,\Phi(\eta)\right].
\end{equation}
Equation \eref{chief} is equivalent to Eq.~\eref{main}, we have not made any approximation. 
Turning to the $\xi\to -\infty$ limit we note that in this regime $\Phi\to 0$ and hence we can expand 
the exponent on the right-hand side of equation \eref{chief}. Keeping only two terms we simplify equation \eref{chief} to
\begin{equation}
\label{chief_2}
\Phi(\xi-v) = \int_{-\infty}^\xi d\eta\,\Phi(\eta).
\end{equation}
{}From \eref{chief_2}, or from equation $\Phi'(\xi-v)= \Phi(\xi)$ obtained by differentiating of
Eq.~\eref{chief_2}, we see that the solution has an exponential form 
\begin{equation}
\label{lambda}
\Phi(\xi) = D\, e^{a \xi}
\end{equation}
Plugging \eref{lambda} into \eref{chief_2} we arrive at the dispersion relation
\begin{equation}
\label{disp}
a\, e^{-a v} = 1
\end{equation}
An elementary analysis of this equation indicates that solutions exist only when $v\leq e^{-1}$. We now invoke the selection principle which asserts that the extremal value, $v=e^{-1}$ in our case, is realized. 

Traveling wave solutions have been investigated in the context of partial differential equations. A few  partial differential equations admitting traveling wave solutions have been deeply studied. One such equation is the celebrated Fisher-KPP equation \cite{Fisher,KPP} for which the selection principle had been proven \cite{KPP,front1} for sufficiently steep initial conditions. (For more recent work see e.g. \cite{bd,evs, evs2}.  A very comprehensive review of traveling wave solutions of non-linear partial differential equations has been given by van Saarloos \cite{van}, a lighter expositions appear in books \cite{mu,B,book}.) More recently, traveling wave solutions have been investigated in the context of nonlinear recurrences arising in the analyses of binary search algorithms \cite{km,bkm,mk}, kinetic theory \cite{van_zon}, and other problems \cite{mk1,mk2,SM,Raissa}; see \cite{mk_rev1,mk_rev2} for a review of the applications of traveling wave techniques to recurrences.  

Asymptotically, the wave front advances at a constant velocity $v=e^{-1}$. The approach to this asymptotic value is rather slow, namely there is a $n^{-1}$ correction in the leading order, resulting in a logarithmic correction to the front position. This correction was first established by Bramson \cite{front1}  for the Fisher-KPP equation; it was subsequently generalized \cite{bd,evs, evs2,van} to more general partial differential equations and to recurrences \cite{km,bkm,mk,mk1,mk2}. This correction generally has the form $\frac{3}{2a}\ln n$. For the selected velocity $v=e^{-1}$, the decay amplitude $a=e$ is implied by dispersion relation \eref{disp}. Taking into account this logarithmic correction we get 
\begin{equation}
\label{front_log}
x_f = e^{-1} n  + \frac{3}{2e}\,\ln n + \mathcal{O}(1)
\end{equation}

It was convenient to think about $x$ and $n$ as space and time coordinates, so that the front of the traveling wave was advancing and we determined $x_f=x_f(n)$. In the original problem, the parameter $x$ is fixed and we are interested in the height $H(x)$ of the cascade tree. The height is essentially the inverse to $x_f=x_f(n)$ which is taken when $x_f=x$. Thus
\begin{equation}
\label{height_log}
H(x) = ex  - \frac{3}{2}\,\ln x + \mathcal{O}(1)
\end{equation}
The height is of course a random quantity. Equation \eref{height_log} gives the average height. In the $x\to\infty$ limit, the average provides a faithful description as it is a growing quantity while the variance remains finite. We haven't proved this assertion, but at least on the physical level of rigor it is obvious: The probability distribution $Pn(x)$ has asymptotically a traveling wave shape with the width of the front remaining finite, and this is essentially equivalent to the finite width of the height distribution.

\section{Size of the Cascade Tree} 
\label{size_tree}

The size $S(x)$ of the cascade tree, that is, the total number of vertices in the tree, is a random variable. Let us compute the average size $\langle S(x)\rangle$. {}From the definition of the continuum cascade model we deduce 
\begin{eqnarray}
\nonumber 
\langle S(x)\rangle &=&  1 + 
\sum _{k=1}^{\infty} \frac{x^k }{k!}\,e^{-x}  \int _0^x \cdots  
\int _0^x  \frac{dy_1}{x}\cdots \frac{dy_k}{x}\left[\langle S(y_1)\rangle + \cdots +\langle S(y_k)\rangle\right]
\\\nonumber
&=&  1 + \sum _{k=1}^{\infty} \frac{x^k }{k!}\,e^{-x}  k\int _0^x \frac{dy}{x}\,\langle S(y)\rangle
\\
&=& 1 + \int _0^x dy\,\langle S(y)\rangle
\label{eq:size}
\end{eqnarray} 
Differentiating \eref{eq:size} we obtain $\frac{d\langle S\rangle}{dx} = \langle S\rangle$, from which 
\begin{equation}
\label{size_1}
\langle S(x)\rangle = e^x
\end{equation}

A similar line of reasoning leads to an integral equation for the second moment
\begin{eqnarray*}
\langle S^2(x)\rangle &=&  1   
+ \sum _{k=1}^{\infty} \frac{x^k }{k!}\,e^{-x}  k\int _0^x \frac{dy}{x}\,\langle S^2(y)\rangle
+ 2\sum _{k=1}^{\infty} \frac{x^k }{k!}\,e^{-x}  k\int _0^x \frac{dy}{x}\,\langle S(y)\rangle\\
&+& \sum _{k=1}^{\infty} \frac{x^k }{k!}\,e^{-x}\,  
\frac{k(k-1)}{2}\int _0^x \frac{dy_1}{x}\,\langle S(y_1)\rangle \int _0^x \frac{dy_2}{x}\,\langle S(y_2)\rangle
\end{eqnarray*}
Using \eref{size_1} we simplify above integral equation to
\begin{equation*}
\langle S^2(x)\rangle = 1 + \int _0^x dy\,\langle S(y)\rangle + 2(e^x-1)+(e^x-1)^2
\end{equation*}
which is solved to yield
\begin{equation}
\label{size_2}
\langle S^2(x)\rangle = 2e^{2x} - e^x
\end{equation}

One can continue and compute 
\begin{equation}
\label{size_3}
\langle S^3(x)\rangle = \frac{15}{4}\,e^{3x} - \frac{11}{4}\,e^x - \frac{3}{2}\,x e^x
\end{equation}
and a few higher moments $\langle S^p(x)\rangle$, but results quickly become very cumbersome. The explicit results \eref{size_1}--\eref{size_3} show that, in contrast to the height, the size of the cascade tree is the random quantity whose limiting distribution (in the $x\to\infty$ limit) remains broad. More precisely, 
in the scaling limit
\begin{equation*}
x\to\infty, \quad S\to\infty, \quad \sigma = e^{-x}S = {\rm finite}
\end{equation*}
the size distribution becomes 
\begin{equation*}
{\rm Prob}[S(x)=S] = e^{-x} F(\sigma)
\end{equation*}
with the limiting distribution being different from the delta function, $F(\sigma)\ne \delta(\sigma - 1)$. 
The normalization requirement together with \eref{size_1}--\eref{size_3} and similar equations for higher moments $\langle S^p(x)\rangle$ show that the moments $M_p=\int_0^\infty d\sigma\, \sigma^p F(\sigma)$ of the limiting distribution are 
\begin{equation*}
M_0 = M_1=1, \quad M_2 = 2, \quad M_3 = \frac{15}{4}\,, \quad M_4 = \frac{34}{3}\,, \quad M_5 = 25,  
\end{equation*}
etc. 

\section{Summary}

We proposed a minimalist model of infinite directed random graphs. The model is a continuum version of a model of finite directed random graphs, known as the cascade model, which has been investigated in the context of food webs and parallel computation.

Our model presumes a total order on the set of vertices. We chose the simplest such set, an interval of length $x$. From each $y\in [0,x]$, directed links to points $y'>y$ are drawn at random according to the Poisson distribution, that is, the points $y'\in (y,x]$ are chosen independently from each other and uniformly with density one. The analysis of this continuum cascade model is actually simpler than the analysis of the discrete cascade model. This is demonstrated by studying the distribution of the length of the longest directed paths starting at the origin (equivalently, the height of the cascade tree with the root at the origin). We employed traveling wave techniques to extract the asymptotic behavior of the length of the longest directed paths in the $x\to\infty$ limit. It will be interesting to understand the limiting distribution of the size of the cascade tree with the root at the origin as well as other properties of the continuum cascade model.

\bigskip
\noindent 
{\bf Acknowledgment} YI  thanks  Joel E. Cohen for helpful comments and discussions.  YI  is  supported in part by US National Science Foundation Grant  DMS 0443803 to Rockefeller University and by JSPS Grant-in-aid for Scientific Research 23540177.


\section*{References}
\begin{thebibliography}{99} 

\bibitem{sr}     
     Solomonoff R and Rapaport A 1959
     Connectivity of random nets
     {\em Bull. Math. Biophys.} {\bf 13} 107--117 

\bibitem{er} 
     Erd\H{o}s P and R\'enyi  A  1960
     On the evolution of random graphs
     {\em Publ. Math. Inst. Hungar. Acad. Sci.} {\bf 5} 17--61

\bibitem{bol} 
    Bollob\'as  B 1985  {\it  Random Graphs}
    (London: Academic Press)

\bibitem{s}
    Shepp L A  1989
    Connectedness of certain random graphs
    {\em Israel J. Math.} {\bf 67} 23--33

\bibitem{dk} Durrett R and Kesten H  1990 
    The critical parameter of the connectedness of some random graphs,   
    A Tribute to Paul Erd\H{o}s,  161--176 (Cambridge: Cambridge University Press)
   
\bibitem{jklp}   
     Janson S, Knuth D E, \L uczak T and Pittel B 1993
     The birth of the giant component
     {\em Rand. Struct. Alg.} {\bf 3}  233--358 

\bibitem{jlr} 
     Janson S, Luczak T  and  Rucinski A  2000 
     {\it  Random Graphs}  (New York: John Wiley \& Sons) 

\bibitem{pf}     
     Flory P F 1953  {\em Principles of Polymer Chemistry}
     (Ithaca: Cornell University Press)

\bibitem{book}  
     Krapivsky P L, Redner S and Ben-Naim E  2010
     {\it  A Kinetic View of Statistical Physics} (Cambridge: Cambridge University Press)

\bibitem{mej}    
      Newman M E J 2002 
      Spread of epidemic disease on networks 
      {\em Phys. Rev. E} {\bf 66}  016128 

\bibitem{gc} 
       Caldarelli G 2007 
       {\em Scale-Free Networks: Complex Webs in Nature and Technology}
       (Oxford: Oxford University Press)

\bibitem{DM03} 
       Dorogovtsev S N and Mendes J F F 2003 
       {\it Evolution of Networks: From Biological Nets to the Internet and WWW} 
       (Oxford: Oxford University Press) 

\bibitem{mejn} 
       Newman M E J 2010
       {\em  Networks: An Introduction} (Oxford: Oxford University Press)

\bibitem{broder}
      Broder A, Kumar R, Maghoul F,  Raghavan P,  Rajagopalan S, Stata R,  Tomkins A and Wiener J 2000
      Graph structure in the web
      {\em Computer Networks} {\bf 33} 309--320

\bibitem{krr}
      Krapivsky P L, Rodgers G J and Redner S   2001
      Degree distributions of growing networks 
      {\em Phys.\ Rev.\ Lett.} {\bf 86} 5401--5404 

\bibitem{kr}
       Krapivsky P L and Redner S   2002
       A statistical physics perspective on web growth
       {\em Computer Networks} {\bf 39} 261--276

\bibitem{c}
     Cohen J  E  1990  A stochastic theory of community food webs VI. 
     Heterogeneous alternatives to the cascade model
     {\em Theor. Popul. Biol.} {\bf 37} 55--90

\bibitem{cbn}
     Cohen J E, Briand F  and Newman C M  1990
     {\it Community food webs: Data and Theory} (New York: Springer-Verlag)

\bibitem{cn}
     Cohen J E and Newman C M 1985 
     A stochastic theory of community food webs: I.  Models and aggregated data
     {\em Proc. R. Soc.  (London)  B} {\bf 224} 421--448

\bibitem{cbn2}
     Cohen J E, Briand F  and Newman C M  1986
     A stochastic theory of community food webs III. Predicted and observed lengths of food chains 
     {\em Proc. R. Soc.  (London)  B} {\bf  228} 317--353

\bibitem{cn2}
    Cohen J E and Newman C M  1986 
    A stochastic theory of community food webs IV. Theory of food chain length in large webs 
    {\em Proc. R. Soc.  (London)  B} {\bf  228} 355--377

\bibitem{n}
   Newman C M  1992
   Chain lengths in certain random directed graphs 
   {\em Rand. Struct. Alg.} {\bf 3} 243--253
 
\bibitem{gel} 
    Gelenbe E,   Nelson R,   Philips T   and Tantawi A 1986 
    An approximation of the processing time for a random graph model of parallel computation, 
    ACM 86 Proceedings of 1986 ACM Fall joint computer conference  
    IEEE Computer Society Press, Los Alamos, 691--697

\bibitem{newman_94}
    Isopi M and Newman C M  1992
    Speed of parallel-processing for random task graphs
    {\em Commun. Pure Appl. Math.} {\bf 47} 361--376

\bibitem{Drmota}
     Drmota M 2009
     {\it Random Trees: An Interplay between Combinatorics and Probability}
     (Wien: Springer)

\bibitem{si} 
    Sibuya M and Itoh Y  1987
    Random sequential bisection and its associated binary tree 
    {\em Ann. Inst. Stat. Math.} {\bf 39} 69--84

\bibitem{ho} 
    Hattori T and Ochiai H  2006  
    Scaling limit of successive approximations for $w'=-w^2$ 
    {\em Funkcialaj Ekvacioj} {\bf 39} 291--319

\bibitem{jn} 
    Janson S and Neininger R  2008 
    The size of random fragmentation trees 
    {\em Probab. Theory  Rel. Fields} {\bf 142} 399--442
 
\bibitem{di} 
    Dutour  Sikiric M  and Itoh Y  2011
    {\it  Random Sequential Packing of Cubes} (London: World Scientific)

\bibitem{r} 
    Robson J M  1979 
    The height of binary search trees, 
    {\em Australian Comput. J.} {\bf 11} 151--153

\bibitem{fo} 
    Flajolet P  and Odlyzko A 1982 
    The average height of binary tree and other simple tree 
    {\em J. Comput. Syst. Sci.}  {\bf 25} 171--213 

\bibitem{d} 
    Devroye L  1986  
    A note on the height of binary search trees 
    {\em J. ACM} {\bf 33} 489--498

\bibitem{km}
     Krapivsky P L  and Majumdar S N  2000
     Traveling waves, front selection, and exact nontrivial exponents
     in a random fragmentation problem 
     {\em Phys.\ Rev.\ Lett.} {\bf 85} 5492--5495

\bibitem{bkm}
     Ben-Naim E, Krapivsky P L  and Majumdar S N  2001
     Extremal properties of random trees 
     {\em Phys.\ Rev.\ E} {\bf 64} 035101

\bibitem{mk}
    Majumdar S N and Krapivsky P L  2002
    Extreme value statistics and traveling fronts: Application to computer science
    {\em Phys.\ Rev.\ E} {\bf 65} 036127
       
\bibitem{sz}
     Szpankowski W  2001 
     {\it Average case analysis of algorithms on sequences} (New York: Wiley)

\bibitem{itoh}
    Itoh Y  2011   Random sequential generation  of intervals 
    for the cascade model of food webs 
    {\em arXiv:1106.4701}
    
\bibitem{Fisher}  
    Fisher R A  1937   The wave of advance of advantageous genes
    {\em Ann. Eugenics} {\bf 7} 355--369

\bibitem{KPP} 
    Kolmogorov A,  Petrovsky I,  and Piskunov N   1937 Mosc. Univ. Bull. Math. A {\bf 1} 1; 
    translated and reprinted in P. Pelce  1988 {\it Dynamics of Curved Fronts} (San Diego: Academic)

\bibitem{front1} 
    Bramson M  1983
    {\it Convergence of Solutions of the Kolmogorov Equation to Traveling Waves} 
    American Mathematical Society, Providence, R.I. 

\bibitem{bd}
     Brunet E and Derrida B  1997  Shift in the velocity of a front due to a cut-off  
     {\em Phys. Rev. E} {\bf 56}  2597--2604.

\bibitem{evs}
     Ebert U and van Saarloos W 1998 
     Universal algebraic relaxation of fronts propagating into an unstable state and implications 
     for moving boundary approximations
     {\em Phys.\ Rev.\ Lett.} {\bf 80} 1650--1653
    
\bibitem{evs2}
      Ebert U and van Saarloos W  2000
      Front propagation into unstable states: universal algebraic convergence towards uniformly 
      translating  pulled fronts 
      {\em Physica D} {\bf 146} 1--99

\bibitem{van}
     van~Saarloos W 2003  
     Front propagation into unstable states 
     {\em Phys. Rep.} {\bf 386} 29--222 

\bibitem{mu}
     Murray J D 1989 
     {\it Mathematical Biology} (New York: Springer-Verlag)

\bibitem{B} 
     Barenblatt G I  1995  
     {\it Scaling, Self-Similarity, and Intermediate Asymptotics}
     (Cambridge: Cambridge University Press)

\bibitem{van_zon}
     van Zon R, van Beijeren H and Dellago Ch 1998
     Largest Lyapunov exponent for many particle systems at low densities
     {\em Phys. Rev. Lett.} {\bf 80} 2035

\bibitem{mk1}
     Majumdar S N and Krapivsky P L  2000
     Extremal paths on a random Cayley tree
     {\em Phys.\ Rev.\ E} {\bf 62} 7735--7742

\bibitem{mk2}
     Majumdar S N and Krapivsky P L  2001
     The dynamics of efficiency: A simple model
     {\em Phys.\ Rev.\ E} {\bf 63} 045101(R) 
       
\bibitem{SM}
     Majumdar S N    
     Traveling front solutions to directed diffusion limited aggregation, digital search trees and 
     the Lempel-Ziv data compression algorithm
     {\em Phys. Rev. E} {\bf 68} 026103    
       
\bibitem{Raissa}
      D'Souza R M, Krapivsky P L and Moore C 2007
      The power of choice in growing trees
      {\em Eur. Phys. J. B} {\bf 59} 535--543

\bibitem{mk_rev1}
     Majumdar S N and Krapivsky P L  2003
     Extreme value statistics and traveling fronts: Various applications
     {\em Physica A} {\bf 318}, 161--170
     
\bibitem{mk_rev2}
     Majumdar S N, Dean D S and Krapivsky P L  2005
     Understanding search trees via statistical physics
     {\em Pramana J. Phys.} {\bf 64}, 1175--1189

\end{thebibliography}
\end{document}